\documentclass{article}
\usepackage{amssymb,latexsym,amscd}

\newcommand{\se}[1]{{\section{#1}} {\setcounter{equation}{0}}}
\newtheorem{theorem}{Theorem}[section]

\newtheorem{prop}{Proposition}[section]

\def\k{{K\"{a}hler }}
\def\ke{{K\"{a}hler-Einstein }}

\def\wp{{Weil-Peterson }}
\begin{document}
\hbadness=10000
\title{{\bf Degeneration of \ke Manifolds I:}\\
{\Large {\bf The Normal Crossing Case}}}
\author{Wei-Dong Ruan\\
Department of Mathematics\\
University of Illinois at Chicago\\
Chicago, IL 60607\\}
\date{October 2002}
\footnotetext{Partially supported by NSF Grant DMS-9703870 and DMS-0104150.}
\maketitle
\begin{abstract}
In this paper we prove that the \ke metrics for a degeneration family of \k manifolds with ample canonical bundles Gromov-Hausdorff converge to the complete \ke metric on the smooth part of the central fiber when the central fiber has only normal crossing singularities inside smooth total space. We also prove the incompleteness of the Weil-Peterson metric in this case. 
\end{abstract}
\se{Introduction}
According to the conjecture of Calabi, on a complex manifold $X$ with ample canonical bundle $K_X$, there should exist a \ke metric $g$. Namely, a metric satisfying ${\rm Ric}_g = -\omega_g$, where $\omega_g$ is the \k form of the \k metric $g$. The existence of such metrics was proved by Aubin and Yau (\cite{Yau2}) using complex Monge-Amp\`{e}re equation. This important result has many applications in \k geometry. Since \ke metric is canonical for the manifold, one would expect its structure to be closely related to the topology and complex geometry of the manifold. Starting with this important result, Yau initiated the program of application of \ke metrics to algebraic geometry (\cite{yau}). It was realized by him the need to study such metrics for quasi-projective manifolds (\cite{yau1}) and their degenerations. The original proof (\cite{Yau2}) was a purely existence result. Later the existence of $C_1(X)<0$ \ke metrics was generalized to complete complex manifolds by Cheng and Yau (\cite{CY}), where other than existence the proof also exhibits the asymptotic behavior of the \ke metric near the infinity boundary.\\

From algebraic geometry point of view, when discussing the compactification of the moduli space of complex manifold $X$ with ample canonical bundle $K_X$, it is necessary to consider a holomorphic degeneration family $\pi: {\cal X} \rightarrow B$ such that $X_t = \pi^{-1}(t)$ are smooth except for $t=0$, and such that the canonical bundle of $X_t$ for $t\not=0$ as well as the dualizing sheaf of $X_0$ are ample. One important case is when ${\cal X}$ is smooth and the central fibre $X_0$ is the union of smooth divisors $D_1,\cdots, D_l$ with normal crossings. As is well known, $K_{\cal X}|_{X_t} = K_{X_t}$ for $t\not=0$ and $K_{X_0}:=K_{\cal X}|_{X_0}$ is the dualizing sheaf of $X_0$. $K_{\cal X}|_{D_i} = K_{D_i} + \sum_{j\not=i}D_j$. The condition that the dualizing sheaf of $X_0$ is ample is equivalent to each line bundle $K_{D_i} + \sum_{j\not=i}D_j$ being ample on $D_i$ for $1\leq i\leq l$.\\

More generally for any index set $I = \{i_0,\cdots, i_k\}$, $K_{\cal X}|_{D_I} = K_{D_I} + \sum_{j\not\in I}D_j$. Recall that for $I\subset J$, there is a natural multi-residue map:\\
\[
r_{IJ}:\ K_{D_I} + \sum_{k\not\in I}D_k \rightarrow K_{D_J} + \sum_{k\not\in J}D_k.
\]\\
A section $\Omega$ of $K_{X_0}$ can be seen as a collection $\{\Omega_I\}$ of sections $\Omega_I$ of $K_{D_I} + \sum_{k\not\in I}D_k$ that satisfy $r_{IJ}(\Omega_I) = \Omega_J$ for $I\subset J$.\\

For $t\not=0$, the works of Aubin and Yau imply the existence of unique \ke metric $g_t$ on $X_t$. For $t=0$, the work of Cheng and Yau implies the existence of a complete \ke metric $g_{0,i}$ on each $D_i\setminus {\rm Sing}(X_0)$. It is interesting to understand the relation between $g_t$ for $t\not=0$ and the collection of complete $g_{0,i}$'s. Such understanding will provide structure results on the \ke metric $g_t$ (which was only known to exist previously) based on the structure of $g_{0,i}$'s. Extending Yau's program, G. Tian made the first important contribution (in \cite{Tian1}) in this direction. He proved (in \cite{Tian1}) that the \ke metric $g_t$ on $X_t$ is Gromov-Hausdorff convergent to the complete \ke metric $g_0 = \{g_{0,i}\}_{i=1}^l$ on the smooth part of $X_0$ under the following three assumptions: \\

(1). the total space ${\cal X}$ is smooth;\\
(2). the central fibre $\displaystyle X_0 = \bigcup_{i=1}^l D_i$ has only normal crossing singularities; \\
(3). any three of the $D_i$'s have empty intersection.\\

Assumption (2) (normal crossing condition) is rather natural. But (1) and (3) are technical requirements that most degeneration families (beyond curves) will not satisfy, and it is desirable to remove them. Among (1) and (3), (1) (the smoothness of the total space ${\cal X}$, or more generally, the base point free condition of the family) is much harder to remove.\\

Tian's proof in \cite{Tian1} can be naturally separated into three steps. The first step is the construction of certain smooth family of background \k metrics $\hat{g}_t$ on $X_t$ and their \k potential volume forms $\hat{V}_t$. (Assumptions (1) and (2) are needed in this step.) The second step is to construct a smooth family of approximate \k metrics $g_t$ with \k form $\displaystyle\omega_t = \frac{i}{2\pi}\partial\bar{\partial}\log V_t$, where $V_t = h\hat{V}_t$ ($h$ is a function on ${\cal X}$) satisfies certain uniform estimate independent of $t$. (Tian needed  assumptions (1), (2) and (3) in this step.) The third step is to start with the smooth family of approximate \k metric $g_t$, using Monge-Amp\`{e}re estimate of Aubin and Yau to derive uniform estimate (independent of $t$) for the smooth family of \ke metrics $g_{E,t}$, which is enough to ensure the Gromov-Hausdorff convergence of the family to the unique complete \ke metric $g_{E,0} = \{g_{0,i}\}_{i=1}^l$ on the smooth part of $X_0$. (No retriction is needed for this step.) The most crucial step is the second step.\\

In this work, we will generalize Tian's result by removing assumption (3). Our proof will follow the general frame work of Tian's and proceed in three steps. The major differences of our approach from Tian's are mainly in the first two steps. In the first step, we provide a very simple construction of background \k metrics $\hat{g}_t$ that is valid even for the more general case (without all three assumptions of Tian as long as the dualizing sheaf of the central fibre is ample). In the second step, we use a construction of $h$ different from Tian's which enables us to get the necessary estimates to carry out the third step without the restriction that any three of the $D_i$'s have empty intersection. (We still need assumptions(1) and (2) here.) Therefore, besides proving our result, we are also able to reduce the program for general case to the construction of function $h$ in the second step. Our first main theorem is the following.\\

\begin{theorem}
\label{aa}
Let $\pi: {\cal X} \rightarrow B$ be a degeneration of \ke manifolds $\{X_t,g_{E,t}\}$ with ${\rm Ric}(g_{E,t}) = - g_{E,t}$. Assume that the total space ${\cal X}$ is smooth and the central fibre $X_0$ is the union of smooth normal crossing hypersurfaces in ${\cal X}$ with ample dualizing line bundle $K_{X_0}$. Then the \ke metrics $g_{E,t}$ on $X_t$ converge to a complete Cheng-Yau \ke metric $g_{E,0}$ in the sense of Cheeger-Gromov.\\
\end{theorem}

In the work of Leung and Lu (\cite{Lu}), they gave a proof (very different from ours) to a result that implies theorem \ref{aa}. We believe our method here can be easily adopted to provide a new (hopefully easier) proof to their result. They followed the general frame work of Tian too. Their first step is the same as Tian's construction of the background \k metrics using cut and paste. In the second step, they constructed $h$ (different from both Tian's and ours) that does not quite satisfy the estimates necessary to carry out the third step as in Tian's work. They compensated this by using a more involved Monge-Amp\`{e}re estimate of Aubin and Yau in the third step. Therefore, although following the same general frame work of Tian and achieving similar result, our approaches in the three steps are all somewhat different from those in \cite{Lu}. In addition, the more precise estimate in our second step enable us to work out the estimate of the \wp metric near degeneration that generalizes Tian's estimate, which implies that the \wp metric is incomplete.\\ 

\begin{theorem}
\label{ab}
The restriction of the \wp metric on the moduli space of complex structures to the degeneration $\pi: {\cal X} \rightarrow B$ is bounded from above by a constant multiple of $\displaystyle\frac{dt\wedge d\bar{t}}{|\log |t||^3|t|^2}$. In particular, \wp metric is incomplete at $t=0$.\\
\end{theorem}

The three steps of our main construction are carried out in sections 2,3,4. Theorem \ref{aa} is proved in section 4. The incompleteness of the \wp metric is discussed in section 5.\\

{\bf Note of notation:} We say $A\sim B$ if there exist constants $C_1,C_2>0$ such that $C_1B \leq A \leq C_2B$.\\

\se{Construction of the background metric}
For construction in this section to work, it is necessary to assume that the dualizing line bundle $K_{\cal X}$ of the total space ${\cal X}$ exists and is ample, which is valid in our situation. Recall that $K_{{\cal X}/B} = K_{\cal X}\otimes K_B^{-1}$ and $K_{X_t}= K_{{\cal X}/B}|_{X_t} \cong K_{\cal X}|_{X_t}$. (The last equivalence is not canonical, depending on the trivialization $K_B \cong {\cal O}_B$. We will use $dt$ to fix the trivialization of $K_B$.) Since $K_{X_t}$ is ample for all $t$, certain multiple $K^m_{X_t}$ will be very ample for all $t$. Equivalently, $K_{\cal X}^m$ is very ample on ${\cal X}$. It is not hard to find sections $\{\Omega_k\}_{k=0}^{N_m}$ of $K_{\cal X}^m$ that determine an embedding $e: {\cal X} \rightarrow {\bf CP}^{N_m}$, such that $\{\Omega_{t,k}\}_{k=0}^{N_m}$ forms a basis of $H^0(K_{X_t})$ for all $t$, where $\Omega_{t,k} = (\Omega_k\otimes (dt)^{-m})|_{X_t}$. $\{\Omega_{t,k}\}_{k=0}^{N_m}$ will determine a family of embedding $e_t: X_t \rightarrow {\bf CP}^{N_m}$ such that $e_t = e|_{X_t}$. Choose the Fubini-Study metric $\omega_{FS}$ on ${\bf CP}^{N_m}$, and define\\
\[
\hat{\omega} = \frac{1}{m}e^*\omega_{FS},\ \ \hat{\omega}_t = \hat{\omega}|_{X_t} = \frac{1}{m}e_t^*\omega_{FS}.
\]\\
Since $K_{\cal X}^m$ is very ample on ${\cal X}$, $\hat{\omega}$ is a smooth metric on ${\cal X}$. The \k potential of $\hat{\omega}$ and $\hat{\omega}_t$ are the volume forms\\
\[
\hat{V} = \left(\sum_{k=0}^{N_m} \Omega_{k}\otimes\bar{\Omega}_{k}\right)^{\frac{1}{m}},\ {\rm and}\ \hat{V}_t = \left(\sum_{k=0}^{N_m} \Omega_{t,k}\otimes\bar{\Omega}_{t,k}\right)^{\frac{1}{m}} = \left.\hat{V} \otimes (dt\otimes d\bar{t})^{-1}\right|_{X_t}.
\]

Choose a tubular neighborhood $U_i$ of $D_i$ and for each $(k+1)$-tuple $I = (i_0,\cdots, i_k)$, denote $U_I = U_{i_1\cdots i_k} = U_{i_1}\cap\cdots\cap U_{i_k}$, which is a tubular neighborhood of $D_I = D_{i_1\cdots i_k} = D_{i_1}\cap\cdots\cap D_{i_k}$. Let $\displaystyle U_I^0 = U_I \setminus \left(\bigcup_{I \subset J} \tilde{U}_J\right)$, where $\tilde{U}_J$ is the closure of slightly shrinked $U_J$. On $U_I^0$, we may choose local coordinate $z = (z_0,\cdots,z_n)$ such that $\displaystyle \prod_{j=0}^kz_j=t$ and $D_{i_j}$ is defined by $z_j=0$ for $0\leq j \leq k$. Since $K_{\cal X}^m$ is ample and therefore base point free, $\hat{V}$ is a non-degenerate smooth volume form on ${\cal X}$. Locally 

\[
\hat{V} = \rho(z)\prod_{j=0}^n dz_jd\bar{z}_j,
\]

where $\rho(z) \sim 1$ is a smooth positive function. Use $(z_1,\cdots,z_n)$ as coordinate on $X_t$, through straightforward computation, we have

\[
\hat{V}_t = \rho(z)\left(\prod_{j=1}^k \frac{dz_jd\bar{z}_j}{|z_j|^2}\right)\left(\prod_{j=k+1}^n dz_jd\bar{z}_j\right).
\]

\se{Construction of the approximate metric}
Fix Hermitian metrics $\|\cdot\|_i$ of the line bundles ${\cal O}(D_i)$ on ${\cal X}$. Let $s_i$ be the section of line bundle ${\cal O}(D_i)$ defining $D_i$. One may choose a suitable trivialization of ${\cal O}(X_0)$ such that $s_1\cdots s_l =t$ on $\cal{X}$. Consider the globally defined functions\\
\[
a_k = \log\|s_k\|_k^2,\ a^2 = \sum_{k=1}^l a_k^2,\ h= \frac{(\log |t|^2)^2}{\displaystyle\prod_{k=1}^la_k^2}.
\]\\
Let $V = h\hat{V}$, then\\
\[
\omega = \frac{i}{2\pi}\partial\bar{\partial}\log V = \hat{\omega} + \frac{i}{2\pi}\partial\bar{\partial}\log h
\]
\[
= \hat{\omega} + \sum_{k=1}^l\frac{2}{a_k}{\rm Ric}(\|\cdot\|_k) + \alpha = \tilde{\omega} + \alpha,
\]

where

\[
\alpha = \frac{i}{\pi}\sum_{k=1}^l\frac{1}{a_k^2}\partial a_k\bar{\partial}a_k,
\]

is always semi-positive. When restricted to $X_t$, we have

\[
\omega_t =  \tilde{\omega}_t + \alpha_t,\ {\rm where}\ \omega_t = \omega|_{X_t},\ \tilde{\omega}_t = \tilde{\omega}|_{X_t},\ \alpha_t = \alpha|_{X_t}.
\]
By possible multiplication by a constant, we may assume $\|s_k\|\leq \delta$ to be small. Then $|a_k|$ will be large and easily $\frac{1}{2}\hat{\omega} \leq \tilde{\omega} \leq 2\hat{\omega}$. Notice that $V_t = h\hat{V}_t$ is the \k potential of $\omega_t$. Assume

\[
e^{-\phi_t} = \frac{\omega_t^n}{V_t}.
\]

\begin{prop}
\label{cb}
$|\phi_t|$ is bounded independent of $t$.
\end{prop}
{\bf Proof:} It is sufficient to verify in each $U_I^0$. For $I=\{i_0,\cdots, i_k\}$, on $U_I^0$, we see that the dominating term of $\omega_t^n$ is $\tilde{\omega}_t^{n-k} \wedge \alpha^k$. More precisely, there is a smooth non-zero function $\rho$ on ${\cal X}$ such that\\
\[
\prod_{i=1}^l\|s_i\|_i^2 = \rho|t|^2 \ \ \ {\rm on}\ {\cal X}.
\]\\
When restricted to $X_t$, we have\\
\[
\sum_{i=1}^l \partial a_i = \partial\log \rho\ \ \ {\rm on}\ X_t.
\]\\
Use this relation to substitute $\partial a_{i_0}$ and $\bar{\partial}a_{i_0}$ in $\omega_t^n$, then the only non-bounded terms will be those involving $\partial a_i$, $\bar{\partial}a_i$ for $i\in I\backslash \{i_0\}$. The dominating term of $\omega_t^n$ is the one involving $\displaystyle\prod_{j(\not=i_0)\in I}\partial a_j\bar{\partial}a_j$, namely

\[
\omega_t^n \sim \tilde{\omega}_t^{n-k} \wedge k!\left(\frac{i}{\pi}\right)^k\left(\sum_{i\in I}a_i^2\right)\left(\prod_{i\in I} a_i^2\right)^{-1}\prod_{i(\not=i_0)\in I}\partial a_i\bar{\partial}a_i.
\]

Choose local coordinate $z=(z_1,z_2,\cdots,z_n)$ on $U^0_I$ such that $D_i\cap U^0_I$ is defined by $z_i=0$ for $i\in I$. Then there are smooth non-zero functions $\rho_i$ such that $\|s_i\|_i^2 = \rho_i|z_i|^2$. Then

\[
\partial a_i = \partial\log \rho_i + \frac{dz_i}{z_i}.
\]

Hence

\[
\omega_t^n \sim k!\left(\frac{i}{\pi}\right)^k\left(\sum_{i\in I}a_i^2\right)\left(\prod_{i\in I} a_i^2\right)^{-1}\left(\prod_{i(\not=i_0)\in I}\frac{dz_id\bar{z}_i}{|z_i|^2}\right)\left(\prod_{i\not\in I}dz_id\bar{z}_i\right).
\]

Since $a_i^2$, $(\log |z_i|^2)^2$ are bounded for $i\not\in I$, we have

\[
\sum_{i\in I}a_i^2 \sim a^2 \sim (\log |t|^2)^2.
\]

For $i\in I$, $a_i^2 \sim (\log |z_i|^2)^2$. Notice that on $U_I^0$,

\[
V_t \sim \frac{(\log |t|^2)^2}{\displaystyle\prod_{i=1}^l a_i^2}\prod_{i=1}^n\frac{dz_id\bar{z}_i}{|z_i|^2} \sim \frac{(\log |t|^2)^2}{\displaystyle\prod_{i\in I} a_i^2}\left(\prod_{i(\not=i_0)\in I}\frac{dz_id\bar{z}_i}{|z_i|^2}\right)\left(\prod_{i\not\in I}dz_id\bar{z}_i\right).
\]

Therefore

\[
e^{-\phi_t} = \frac{\omega_t^n}{V_t} \sim \frac{\sum_{i\in I}a_i^2}{(\log |t|^2)^2} \sim 1.
\]

We have $|\phi_t|$ is bounded. 
\begin{flushright} \rule{2.1mm}{2.1mm} \end{flushright}
Let $g_t$ denote the \k metric corresponding to the \k form $\omega_t$, then we have\\
\begin{prop}
\label{cc}
The curvature of $g_t$ and its derivatives are all uniformly bounded with respect to $t$.\\
\end{prop}
{\bf Proof:} On a Riemannian manifold $(M,g)$, we call a basis $\{v_i\}$ proper if the corresponding metric matrix satisfies $C_1(\delta_{ij}) \leq (g_{ij}) \leq C_2(\delta_{ij})$ for $C_1,C_2>0$. To verify that the Riemannian metric $g$ has bounded curvature, it is enough to find a proper basis $\{v_i\}$ such that the second derivatives of $g_{ij}$ and $C^1$ norm of the coefficients of $[v_i,v_j]$ with respect to the basis $\{v_i\}$ are all bounded.\\

For $I = \{0,\cdots,k\}$, in $U_I^0$, we have coordinate $(z_0,\cdots,z_n)$ satisfying $z_i|_{D_i}=0$ for $i\in I$. We will restrict to the part of $U_I^0$ where $a_0^2$ is the largest among all $a_i^2$, and take $z=(z_1,\cdots,z_n)$ as local coordinate for $X_t$. Let $W_i = a_iz_i\displaystyle\frac{\partial}{\partial z_i}$ for $i(\not=0)\in I$ and $W_i = \displaystyle\frac{\partial}{\partial z_i}$ for $i\not\in I$. In $U_I$, the metric can be written as

\[
g_t = \tilde{g}_t + \alpha_t,\ \ {\rm where}\ \alpha_t =  \frac{i}{\pi}\sum_{i=0}^k\frac{1}{a_i^2}\partial a_i\bar{\partial}a_i.
\]

It is straightforward to check that the basis $\{W_i,\bar{W}_i\}_{i=1}^n$ is proper. Namely $C_1(\delta_{ij}) \leq (g_{i\bar{j}}) \leq C_2(\delta_{ij})$ for some $C_1,C_2>0$, where $(g_{i\bar{j}})$ denotes the metric matrix with respect to the basis $\{W_i,\bar{W}_i\}_{i=1}^n$. (For the upper bound estimate, we need $\frac{a_i}{a_0}$ to be bounded, which is due to our restriction to the part of $U_I^0$ where $a_0^2$ is the largest among all $a_i^2$.)\\ 

For $j\in I$, $\|s_j\|^2 = \rho_j|z_j|^2.$

\[
W_i(a_j) = \frac{W_i(\|s_j\|^2)}{\|s_j\|^2} = \frac{W_i(\rho_j)}{\rho_j} + \frac{W_i(|z_j|^2)}{|z_j|^2}.
\]

\[
W_i(a_j) = a_i(z_i\frac{\partial\log \rho_j}{\partial z_i} + \delta_{ij})\ \ {\rm for}\ i\in I.
\]

\[
W_i(a_j) = \frac{\partial\log \rho_j}{\partial z_i}\ \ {\rm for}\ i\not\in I.
\]

The functions

\begin{equation}
\label{ca}
\frac{a_i}{a},\ \frac{\log|t|^2}{a},\ \frac{\log|t|^2}{a_0},\ \frac{1}{a_i},\ a_iz_i,\ a_i\bar{z}_i,\ \frac{a_i}{a_0},\ \ {\rm for}\ i\in I
\end{equation}

are all bunded in the part of $U_I^0$ where $a_0^2$ is the largest among all $a_i^2$. Above computations imply that the derivatives of functions in (\ref{ca}) with respect to $\{W_i,\bar{W}_i\}_{i=1}^n$ will be smooth functions of terms in (\ref{ca}) and other smooth bounded terms. Therefore they are bounded.\\

It is straightforward to check that $g_{i\bar{j}}$ and coefficients of $[W_i,W_j]$, $[W_i,\bar{W}_j]$, $[\bar{W}_i,\bar{W}_j]$ with respect to the basis $\{W_i,\bar{W}_i\}_{i=1}^n$ are all smooth functions of terms in (\ref{ca}) and other bounded smooth terms. Consequently, their any derivatives with respect to $\{W_i,\bar{W}_i\}_{i=1}^n$ are also smooth functions of terms in (\ref{ca}) and other bounded terms, therefore, are all bounded.
\begin{flushright} \rule{2.1mm}{2.1mm} \end{flushright}
\begin{prop}
\label{cd}
For any $k$, $\|\phi_t\|_{C^k,g_t}$ is uniformly bounded with respect to $t$.\\
\end{prop}
{\bf Proof:} Similar as in the proof of previous proposition, for $I = \{0,\cdots,k\}$, in the part of $U_I^0$ where $a_0^2$ is the largest among all $a_i^2$, $\phi_t$ is a smooth function of terms in (\ref{ca}) and other smooth bounded terms. Consequently, $W_i(\phi_t)$ is also a smooth function of terms in (\ref{ca}) and other bounded terms. By induction, all higher derivatives of $\phi_t$ with respect to $\{W_i\}$ will be a smooth function of terms in (\ref{ca}) and other bounded terms. Therefore they are bounded.
\begin{flushright} \rule{2.1mm}{2.1mm} \end{flushright}

\se{Construction of \ke metric via complex Monge-Amp\`{e}re}
In this section, we will use the same notions as in the previous sections. In \cite{Tian1}, using the Monge-Amp\`{e}re estimate of Aubin and Yau, Tian essentially proved the following.\\

\begin{theorem}
\label{db}
(Tian) Assume that $\phi_t$, the curvature of $g_t$ and their multi-derivatives are all bounded uniformly independent of $t$, then the \ke metric $g_{E,t}$ on $X_t$ will converge to the complete Cheng-Yau \ke metric $g_{E,0}$ on $X_0\setminus{\rm Sing}(X_0)$ in the sense of Cheeger-Gromov: there are an exhaustion of compact subsets $F_\beta \subset X_0\setminus{\rm Sing}(X_0)$ and diffeomorphisms $\psi_{\beta,t}$ from $F_\beta$ into $X_t$ satisfying:\\
(1) $\displaystyle X_t\setminus \bigcup_{\beta=1}^\infty \psi_{\beta,t}(F_\beta)$ consists of finite union of submanifolds of real codimension 1;\\
(2) for each fixed $\beta$, $\psi_{\beta,t}^*g_{E,t}$ converge to $g_{E,0}$ on $F_\beta$ in $C^k$-topology on the space of Riemannian metrics as $t$ goes to $0$ for any $k$.\\
\end{theorem}
\begin{flushright} \rule{2.1mm}{2.1mm} \end{flushright}
{\bf Proof of theorem \ref{aa}:} This theorem is a direct corollary of theorem \ref{db} and propositions \ref{cb}, \ref{cc}, \ref{cd}. 
\begin{flushright} \rule{2.1mm}{2.1mm} \end{flushright}
It is easy to see that our construction actually implies the following asymptotic description of the family of \ke metrics.\\

\begin{theorem}
\label{da}
\ke metric $g_{E,t}$ on $X_t$ is uniformly quasi-isometric to the explicit approximate metric $g_t$. More precisely, there exist constants $C_1,C_2>0$ independent of $t$ such that $C_1 g_t \leq g_{E,t}\leq C_2 g_t$.\\
\end{theorem}
{\bf Proof:} 
The uniform $C^0$-estimate of the complex Monge-Amp\`{e}re equations implies that $C_1 \omega_t^n \leq \omega_{E,t}^n\leq C_2 \omega_t^n$ for some $C_1,C_2>0$. The uniform $C^2$-estimate of the complex Monge-Amp\`{e}re equations implis that ${\rm Tr}_{g_t} g_{E,t}$ is uniformly bounded from above. Combining these two estimates, we get our conclusion. 
\begin{flushright} \rule{2.1mm}{2.1mm} \end{flushright}

\se{\wp metric near degeneration}
{\bf Example:} In general, we consider $(D^*)^{n+1}$, where $D^* = \{0<|z|<1\}$, with \ke metric\\
\[
\omega = \sum_{k=0}^n\frac{i}{\pi}\frac{dz_k\wedge d\bar{z_k}}{|z_k|^2(\log|z_k|^2)^2} = \sum_{k=0}^n\frac{i}{\pi}\frac{\partial a_k\wedge \bar{\partial}a_k}{a_k^2},
\]\\
where
\[
a_k = \log|z_k|^2,\ a^2 = \sum_{k=0}^n a_k^2.
\]\\
We are interested in the hypersurface\\
\[
X_t = \{\prod_{k=0}^n z_k =t\} \subset (D^*)^{n+1}.
\]\\
Let $\omega_t = \omega|_{X_t}$, then\\
\[
\omega_t = \sum_{k=1}^n\frac{i}{\pi}\frac{dz_k\wedge d\bar{z_k}}{|z_k|^2(\log|z_k|^2)^2} + \frac{i}{\pi}\frac{1}{\left(\log |t|^2 - \displaystyle\sum_{k=1}^n\log|z_k|^2\right)^2}(\sum_{k=1}^n\frac{dz_k}{z_k})\wedge (\sum_{k=1}^n\frac{d\bar{z_k}}{\bar{z_k}}).
\]\\
Let\\
\[
W= \frac{\nabla \log t}{|\nabla \log t|^2},
\]\\
then $\displaystyle\pi_*W = t\frac{d}{dt}$. It is straightforward to derive that\\
\[
W = \sum_{k=0}^n \frac{a_k^2}{a^2} z^k\frac{\partial}{\partial z_k}.
\]\\
\[
\bar{\partial} W = \sum_{k=0}^n \bar{\partial}\left(\frac{a_k^2}{a^2}\right) z^k\frac{\partial}{\partial z_k}
\]\\
is a natural representative of Kodaira-Spencer deformation class in the Dolbeaut cohomology $H^1(T_X)$.
\begin{flushright} \rule{2.1mm}{2.1mm} \end{flushright}
For $c <1$ consider a smaller domain $(D^*_c)^{n+1}$, where $D^*_c = \{0<|z|<c\}$. For $\omega_t$ and $W$ as in the previous example, we have\\
 
\begin{prop}
\[
\frac{\displaystyle\int_{X_t\cap(D^*_c)^{n+1}}\|\bar{\partial} W\|^2 \omega^n_t}{\displaystyle\int_{X_t\cap(D^*_c)^{n+1}} \omega^n_t} = \frac{1}{|\log |t|^2|^3}2n|\log c^2|(1+\frac{\pi}{2})(1+O(\epsilon|\log \epsilon|)).
\] 
\end{prop}
{\bf Proof:} 
Straightforward computation gives us\\
\[
\omega^n_t = n!\left(\frac{i}{\pi}\right)^n \frac{a^2}{\displaystyle\prod_{k=0}^na_k^2} \prod_{k=1}^n\frac{dz_k\wedge d\bar{z_k}}{|z_k|^2} = n!\left(\frac{1}{\pi}\right)^n \frac{a^2}{\displaystyle\prod_{k=0}^na_k^2} \prod_{k=1}^n(da_k\wedge d\theta_k).
\]\\
We may compute the volume of $X_t\cap(D^*_c)^{n+1}$.\\
\[
\int_{X_t\cap(D^*_c)^{n+1}} \omega^n_t = n!2^n \int_{\tiny{\begin{array}{c}a_k\leq \log c^2,\ 1\leq k \leq n\\ \displaystyle\sum_{k=1}^n a_k \geq \log |t|^2 -\log c^2\end{array}}} \frac{a^2}{\displaystyle\left(\prod_{k=1}^na_k^2\right)\left(\log |t|^2 - \sum_{k=1}^na_k\right)^2} \prod_{k=1}^nda_k
\]
\[
 = n!2^n \frac{1}{|\log |t|^2|^n}\int_{\tiny{\begin{array}{c}\\b_k\geq \epsilon,\ 1\leq k \leq n\\ \displaystyle\sum_{k=1}^n b_k \leq 1-\epsilon\end{array}}} \frac{b^2}{\displaystyle\left(\prod_{k=1}^nb_k^2\right)\left(1 - \sum_{k=1}^nb_k\right)^2} \prod_{k=1}^ndb_k
\]
\[
 = (n+1)!2^n \frac{1}{|\log |t|^2|^n}\int_{\tiny{\begin{array}{c}\\b_k\geq \epsilon,\ 1\leq k \leq n\\ \displaystyle\sum_{k=1}^n b_k \leq 1-\epsilon\end{array}}} \left(\prod_{k=1}^nb_k^2\right)^{-1} \prod_{k=1}^ndb_k
\]\\
where\\
\[
\epsilon = \frac{\log c^2}{\log |t|^2},\ b_k = \frac{a_k}{|\log |t|^2|}.
\]\\
From this, it is straightforward to derive that\\
\[
\frac{(n+1)!2^n}{|\log |t|^2|^n}\prod_{k=1}^n\int_{\epsilon}^{\frac{1}{n}} \frac{1}{b_k^2} db_k \leq \int_{X_t\cap(D^*_c)^{n+1}} \omega^n_t\leq  \frac{(n+1)!2^n}{|\log |t|^2|^n}\prod_{k=1}^n\int_{\epsilon}^1 \frac{1}{b_k^2} db_k
\]\\
\[
\frac{(n+1)!2^n}{|\log c^2|^n}(1-\frac{\epsilon}{n})^n \leq \int_{X_t\cap(D^*_c)^{n+1}} \omega^n_t\leq  \frac{(n+1)!2^n}{|\log c^2|^n}(1-\epsilon)^n.
\]

Notice

\[
\|\bar{\partial} W\|^2 = \sum_{k=0}^n \frac{1}{a_k^2}\left\|\bar{\partial}\left(\frac{a_k^2}{a^2}\right)\right\|^2
\]\\
\[
\bar{\partial}\left(\frac{a_k^2}{a^2}\right) = \frac{2a_k}{a^2}\bar{\partial}a_k - \frac{2a_k^2}{a^4}\sum_{j=0}^na_j\bar{\partial}a_j
\]\\
It is straightforward to compute\\
\[
\left\|\bar{\partial}\left(\frac{a_k^2}{a^2}\right)\right\|^2 = \frac{4a_k^4}{a^8}(a^4 -2a_k^2a^2 + \sum_{j=0}^na_j^4)
\]\\
\[
\|\bar{\partial} W\|^2 = \frac{4}{a^6}(a^4 -\sum_{j=0}^na_j^4) = \frac{4}{a^6}(\sum_{i\not=j}^na_i^2a_j^2).
\]\\
\[
\int_{X_t\cap(D^*_c)^{n+1}}\|\bar{\partial} W\|^2 \omega^n_t = n!2^n \int_{\tiny{\begin{array}{c}a_k\leq \log c^2,\ 1\leq k \leq n\\ \displaystyle\sum_{k=1}^n a_k \geq \log |t|^2 -\log c^2\end{array}}} \frac{4\displaystyle\left(\sum_{i\not=j}^na_i^2a_j^2\right) \prod_{k=1}^nda_k}{a^4\displaystyle\left(\prod_{k=1}^na_k^2\right)\left(\log |t|^2 - \sum_{k=1}^na_k\right)^2}
\]
\[
 = n!2^n \frac{1}{|\log |t|^2|^{n+2}}\int_{\tiny{\begin{array}{c}\\b_k\geq \epsilon,\ 1\leq k \leq n\\ \displaystyle\sum_{k=1}^n b_k \leq 1-\epsilon\end{array}}} \frac{4\displaystyle\left(\sum_{i\not=j}^nb_i^2b_j^2\right)}{b^4\displaystyle\left(\prod_{k=1}^nb_k^2\right)\left(1 - \sum_{k=1}^nb_k\right)^2} \prod_{k=1}^ndb_k
\]
\[
 = (n+1)!2^{n+1}n \frac{1}{|\log |t|^2|^{n+2}}\int_{\tiny{\begin{array}{c}\\b_k\geq \epsilon,\ 1\leq k \leq n\\ \displaystyle\sum_{k=1}^n b_k \leq 1-\epsilon\end{array}}} \frac{1}{b^4\displaystyle\left(\prod_{k=1}^{n-1}b_k^2\right)} \prod_{k=1}^ndb_k
\]
\[
= \frac{(n+1)!2^{n+1}n}{|\log c^2|^{n-1}}\frac{1}{|\log |t|^2|^3} \int_0^1 \frac{db_n}{(b_n^2 + (1-b_n)^2)^2} \prod_{k=1}^{n-1}\int_1^{+\infty} \frac{dx_k}{x_k^2} (1+O(\epsilon|\log \epsilon|)).
\]  
\[
= \frac{1}{|\log |t|^2|^3}\frac{(n+1)!2^{n+1}n}{|\log c^2|^{n-1}}(1+\frac{\pi}{2})(1+O(\epsilon|\log \epsilon|)).
\]\\
\[
\frac{\displaystyle\int_{X_t\cap(D^*_c)^{n+1}}\|\bar{\partial} W\|^2 \omega^n_t}{\displaystyle\int_{X_t\cap(D^*_c)^{n+1}} \omega^n_t} = \frac{1}{|\log |t|^2|^3}2n|\log c^2|(1+\frac{\pi}{2})(1+O(\epsilon|\log \epsilon|)).
\] 
\begin{flushright} \rule{2.1mm}{2.1mm} \end{flushright}
With respect to the approximate \k metric $g$ and parametrizing function $t$ on ${\cal X}$, we can similarly define $W$. $\bar{\partial} W$ similarly represents the Kodaira-Spencer deformation class. We have

\begin{prop}
\label{ea}
There exists a constant $C1,C_2>0$ independent of $t$ such that

\[
\frac{C_1}{|\log |t|^2|^3}\int_{X_t} \omega^n_t \leq \int_{X_t}\|\bar{\partial} W\|_{g_t}^2 \omega^n_t \leq \frac{C_2}{|\log |t|^2|^3}\int_{X_t} \omega^n_t.
\]
\end{prop}
{\bf Proof:} 
Locally in each $U_I^0$, we will use similar coordinate and proper basis $\{W_i,\bar{W}_i\}_{i=0}^n$ as in the proof of proposition \ref{cc}. Then the dual basis is $\{\beta_i,\bar{\beta}_i\}_{i=1}^n$, where $\beta_i = \frac{dz_i}{a_iz_i}$, $a_i = \log |z_i|^2$ for $i\in I$ and $\beta_i = dz_i$ for $i\not\in I$. We have

\[
\omega = \sum_{i\in I} \beta_i\bar{\beta}_i + \sum_{i,j\not\in I} g_{i\bar{j}}dz_id\bar{z}_j + O(a_i^{-1},|z_i|a_i,i\in I).
\]

The term $O(a_i^{-1},|z_i|a_i,i\in I)$ is a smooth function of bounded terms in (\ref{ca}), therefore the multi-derivatives of it with respect to $\{W_i,\bar{W}_i\}_{i=0}^n$ are also of order $O(a_i^{-1},|z_i|a_i,i\in I)$, which are small in $U_I^0$. It is straightforward to derive that

\[
W = \sum_{i\in I} \frac{a_i^2}{a^2} z^i\frac{\partial}{\partial z_i} + \frac{1}{a}O(a_i^{-1},|z_i|a_i,i\in I).
\]

\[
\bar{\partial} W = \sum_{i\in I} \bar{\partial}\left(\frac{a_i^2}{a^2}\right) z^i\frac{\partial}{\partial z_i} + \frac{1}{a}O(a_i^{-1},|z_i|a_i,i\in I).
\]

Apply proposition \ref{ea}, we can find $C_1,C_2>0$ independent of $t$ such that 

\[
\frac{C_1}{|\log |t|^2|^3}\int_{U_I^0\cap X_t} \omega^n_t \leq \int_{U_I^0\cap X_t}\|\bar{\partial} W\|_{g_t}^2 \omega^n_t \leq \frac{C_2}{|\log |t|^2|^3}\int_{U_I^0\cap X_t} \omega^n_t.
\]

Conbine these estimates, we get the statement of the proposition.
\begin{flushright} \rule{2.1mm}{2.1mm} \end{flushright}
{\bf Proof of theorem \ref{ab}:} As pointed out in \cite{Tian1}

\[
g_{WP}\left.\left(\frac{d}{dt},\frac{d}{dt}\right)\right|_{X_t} = \int_{X_t}\left\|H\left(\frac{d}{dt}\right)\right\|^2_{g_{E,t}} \omega_{E,t}^n
\]

where $\displaystyle H\left(\frac{d}{dt}\right)$ denote the harmonic representative of the Kodaira-Spencer deformation class. As mentioned earlier, such class can also be represented by $\displaystyle \frac{\bar{\partial} W}{t}$. Apply proposition \ref{ea} and theorem \ref{da}, we have

\[
\int_{X_t}\left\|H\left(\frac{d}{dt}\right)\right\|^2_{g_{E,t}} \omega_{E,t}^n \leq \int_{X_t}\left\|\frac{\bar{\partial} W}{t}\right\|^2_{g_{E,t}} \omega_{E,t}^n \leq C\int_{X_t}\left\|\frac{\bar{\partial} W}{t}\right\|^2_{g_t} \omega_t^n \leq \frac{C}{|\log |t||^3|t|^2}.
\]
\begin{flushright} \rule{2.1mm}{2.1mm} \end{flushright}
{\bf Remark:} $W$ in this section can be used to simplify the proof of Tian's theorem \ref{db}. Let $\mu_t: X_t \rightarrow X_0$ be the map generated by the inverse flow of $W$. Then the restriction $\mu_t^0: X_t\setminus \mu_t^{-1}({\rm Sing}(X_0)) \rightarrow X_0\setminus{\rm Sing}(X_0)$ is a diffeomorphism. Therefore $\psi_t=(\mu_t^0)^{-1}: X_0\setminus{\rm Sing}(X_0) \rightarrow X_t\setminus \mu_t^{-1}({\rm Sing}(X_0))$ is a diffeomorphism. In particular, $\psi_{\beta,t} := \psi_t|_{F_\beta}$ is a diffeomorphism. With such choice of $\psi_{\beta,t}$, $\displaystyle X_t\setminus \bigcup_{\beta=1}^\infty \psi_{\beta,t}(F_\beta)$ can be naturally identified with $\mu_t^{-1}({\rm Sing}(X_0))$, which is naturally a finite union of submanifolds of real codimension 1. Another nice thing about $\psi_t$ is that it is a symplectomorphism.\\

\ifx\undefined\bysame
\newcommand{\bysame}{\leavevmode\hbox to3em{\hrulefill}\,}
\fi

\end{document}